\documentstyle{amsppt}
\document
\input amstex
\loadbold

\def\f{{\bold f}}
\def\al{{\alpha}}

\def\p{{\frak p}}

\def\ra{{\rightarrow}}

\def\ZZ{{\Bbb Z}}
\def\QQ{{\Bbb Q}}
\def\I{{\bold I}}
\def\Qbar{{\overline{\Bbb Q}}}
\def\RR{{\Bbb R}}

\def\bfP{{\bold{f}_P}}
\def\bftP{{\tilde{\bold{f}}_P}}

\def\bfpP{{\bold{f}_P^{prim}}}
\def\hoI{{\bold{h}^o_{\bold I}}}
\def\Hf{{\bold{H}_{\bold{f}}}}

\def\p1{{\pmatrix p&0\\ 0&1 \endpmatrix}}

\define\isoarrow{{\overset\sim\to\longrightarrow}}
\def\-{{-1}}
%\baselineskip=16pt

%\magnification\magstep1
\NoBlackBoxes

\vskip .7in

\centerline {\bf $p$-adic measures and square roots of
special values of triple product $L$-functions}
\vskip .15in
\centerline{\bf by}
\vskip .15in
\tenpoint
\centerline{\bf Michael Harris\footnote{Institut de Mathematiques de
Jussieu-U.M.R. 7586 du CNRS. Supported in part by the National Science
  Foundation, through Grant DMS-9203142.}}
\centerline{U.F.R. de Math\'ematiques}
\centerline{Universit\'e Paris 7}
\centerline{2, Pl. Jussieu}
\centerline{75251 Paris Cedex 05, FRANCE}
\vskip .15in
\centerline{\bf Jacques Tilouine\footnote{Membre, Institut Universitaire de
France}}
\centerline{Universit\'e Paris-XIII-Institut Galil\'ee}
\centerline{Math\'ematiques - B\^at. B}
\centerline{Avenue Jean Baptiste Cl\'ement}
\centerline{93430 Villetaneuse, FRANCE}
\vskip .15in
\hrule
\vskip .3in

%\baselineskip=16pt

%\magnification\magstep1

\centerline {\bf Introduction}

Let $p$ be a prime number.
In this note, we combine the methods of Hida with the results of [HK1]
to define a $p$-adic analytic function, the squares of whose special
values are related to the
values of triple product $L$-functions at their centers of symmetry.
More precisely, let $f$, $g$, and $h$ be classical normalized cuspidal
Hecke eigenforms
of level $1$ and (even) weights $k$, $\ell$, and $m$, respectively,
with $k \geq \ell \geq m$; assume $k \geq \ell + m$.
Let $L(s,f,g,h)$ be the triple
product $L$-function [G1, G2, PSR]; its center
of symmetry is the point $s = \frac{k + \ell + m - 2}{2}$.  Let
$<\bullet,\bullet>_k$ be the normalized Petersson inner product for
modular forms of weight $k$.  Let $\QQ\{f,g,h\}$ be
the field generated over $\QQ$ by the Fourier coefficients of  $f$, $g$,
and $h$.
Using the integral representation for $L(s,f,g,h)$ [op. cit],
Kudla and one of the authors have shown that the quotient
$$
\frac{L(\frac{k + \ell + m - 2}{2},f,g,h)}{\pi^{2k}<f,f>_k^2 \cdot
C(k,\ell,m) }
$$
is a square in $\QQ\{f,g,h\}$.  Here
$C(k,\ell,m) \in \QQ$ is a universal constant, depending
only on $k$, $\ell$, and $m$.  We construct $p$-adic measures which
interpolate the square root of this quotient, as (the
ordinary eigenform associated to) $f$ varies in a Hida family $\bold{f}$.
These are actually {\it generalized measures}, in the sense of
[H1,II]: elements of the finite normal
algebra extensions of the Iwasawa algebra which arise in Hida's theory of
the ordinary Hecke algebra.
In the simplest case, we obtain the following formula (cf. Theorem 2.2.8):
$$
(\frac{D_H(\bold{f},g,h)(k)}{H(k) \cdot K(k)})^2 =
  \frac{L(\frac{w+1}{2},f_k,g,h)}{\pi^{2k}\cdot <f_k,f_k>_k^2 \cdot
C(k,\ell,m)}. $$
Here $f_k$ is (the primitive form associated to) the specialization in
weight $k$ of $\bold{f}$, $D_H(\bold{f},g,h)$ is the analytic
function associated to our $p$-adic measure, $w = k + \ell + m - 3$, and
$H(k)$ and $K(k)$ are
normalizing factors depending on $\bold{f}$.

The construction of this measure is a modification
of Hida's approach to $p$-adic interpolation of Rankin products [H1,II].
Indeed, when the cusp form $g$ is replaced by an appropriate Eisenstein
series,
Hida's method defines a three-variable Rankin product: one variable for the
value of $s$ and the other two for the $p$-adic variation of $f$ and $h$.
However, there is a subtle difference between our $p$-adic construction and
that of Hida.
Let $Z = \ZZ_p^{\times}$ and let $X$ be a
$p$-adic manifold on which $Z$ acts.  Starting with measures $dE$ on $Z$ and
$\mu$ on $X$, Hida constructed the convolution $d(E*\mu)$ by the formula
$$\int_{Z\times X} \phi(z,x) d(E*\mu) = \int_{Z\times X} \phi(z,z^{-1}x)
dE(z)d\mu(x).$$
The measure $dE$ takes values in the space of $p$-adic modular forms.
The fact that $dE$ is supported on $\ZZ_p^{\times}$ forces
$\int_Z \phi(z) dE(z)$ to be of level divisible by $p$.  Hence
the $p$-adic symmetric square $L$-function constructed in
[H3] comes with an Euler factor at $p$ that gives a trivial zero
at $s = 1$.  In our present construction, we instead form the {\it product}
of a measure and a fixed function.  The calculation leading to (2.2.7)
yields the
correct Euler factor at $p$, and thus to the $p$-adic
interpolation formula (2.2.9).  The same idea was used by one of the
authors,
together with R. Greenberg, to
construct a modified symmetric square $L$-function, with the trivial zero
removed,
and thus to obtain a formula for the derivative at $s = 1$ of Hida's
symmetric
square $L$-function at $s = 1$, with arithmetically interesting
implications (see [GT], [HTU]).

The article [HK1] also considers the central critical values when
$k < \ell + m$.  The $p$-adic interpolation
of the critical values of such triple products, as opposed
to the square roots, has recently been obtained by B\"ocherer and
Panchishkin when each of the rank $2$ motives $M(f)$, $M(g)$, and $M(h)$
associated to $f$, $g$, and $h$ is ordinary, cf. [P].

The notion that the square roots of central critical values of
$L$-functions
should have $p$-adic interpolations seems to have first arisen in
connection with
the thesis of A. Mori [M1, M2].  Mori showed that, if $f$ is a holomorphic
modular form and $K$ is an imaginary quadratic field in which $p$ splits,
then the value of $f$ at Heegner points associated to $K$, suitably
normalized, is naturally an Iwasawa function.  When $f$ is a new form,
it should be possible to use Waldspurger's results in [W] to show that
this Iwasawa function $p$-adically interpolates the square roots of
the central critical values of the $L$-functions $L(s,f_K,\chi)$,
where $f_K$ is the base change of $f$ to $K$, as above,
and $\chi$ runs through a continuous family of algebraic
Hecke characters for $K$.  Other examples have been studied
by Stevens [St], Hida [H5], Sofer [So] and Villegas [V].

The first author would like to take this opportunity to thank
his colleagues at the
Universit\'e de Paris-Sud at Orsay, where most of this article
was written in 1993, for providing a uniquely congenial working
environment.  Both authors would like to thank the referee,
whose careful reading uncovered a significant error in
the previous version, and whose suggestions led to a number of
improvements.  In particular, section 1.4 has been thoroughly
rewritten on the basis of the referee's suggestions.

\bigskip
\medskip

\centerline{\bf 1. p-adic measures associated to three modular forms}
\medskip

\noindent {\bf (1.1)} {\it Review of p-adic modular forms.}
Let $\Cal O$ be an algebra of finite-rank over $\Bbb{Z}_p$,
$\Cal{K} = \Cal{O}\otimes_{\ZZ_p}\QQ_p$.
Let $N$ be a positive integer prime to $p$ and let $Z_N =
\varprojlim (\Bbb{Z}/Np^{\alpha}\Bbb{Z})^{\times}$; $\overline{M} =
\overline{M}(N,\Cal O)$
denotes the complete $\Cal O$-algebra of $\Cal O$-valued ($p$-adic)
modular forms of level $Np^{\infty}$ and any weight (i.e., $p$-adic
modular forms in the sense of Katz; see [H1 I or II, \S 1]).  For any
integers $k, \ell \geq 0$,
and any ring $R$, we let $M_{\ell}(Np^{k},R)$ denote the module of
modular forms of weight $\ell$ for $\Gamma_1(Np^{k})$ -- we
say of level $Np^{k}$, for short --  whose $q$-expansion at
the cusp at infinity lies in $R[[q]]$.  Let
$$M_{\ell}(Np^{\infty},R) = \bigcup_{k} M_{\ell}(Np^{k},R).$$
Similarly, we define $S_{\ell}(Np^{k},R)$, $S_{\ell}(Np^{\infty},R)$,
and $\overline{S} = \overline{S}(N,\Cal O)$ to be the corresponding modules
of cusp forms.  If $g \in \overline{S}$, we write its $q$-expansion
$g = \sum_{n=1}^{\infty} a(n,g) q^n$, with $a(n,g) \in \Cal{O}$ for all
$n$.
We use the same notation for modular forms with complex Fourier
coefficients.
\smallskip

\noindent The following operators on $M_{\ell}(Np^{\infty},\Cal O)$ are
standard,
preserve the submodules of cusp forms, and extend to the completion
$\overline{M}$ (cf. [H3, \S 1] and [Go] for details):

\noindent (1.1.1) For any prime number $\lambda$, the Hecke operators
$T(\lambda)$ and, for $\lambda$ relatively
prime to $N$, $T(\lambda,\lambda)$,
whose action on the
$q$-expansion is given by [H2, (1.13 a)]; more generally, for any
$\lambda$ relatively
prime to $N$, $T(\lambda)$ and $T(\lambda,\lambda)$ can be defined by
the usual formulas.

\noindent (1.1.2)  For any prime number $\lambda$ relatively
prime to $Np$, the {\it diamond operators} $<\lambda>_{\ell}$, whose action
on
$M_{\ell}(Np^{\infty},\Cal O)$ is given by
$$
<\lambda>_{\ell} = \frac{T(\lambda)^2 - T(\lambda^2)}{\lambda^{\ell - 1}}
= \frac{T(\lambda,\lambda)}{\lambda^{\ell - 1}};
$$
we let $<\lambda> = <\lambda>_0$.

\noindent (1.1.3)  By continuity, the map $\lambda\mapsto <\lambda>$ 
extends to an action of the
group
$Z_N$ on $M(Np^{\infty},\Cal O)$ given by $<z>=z_p\cdot T(z,z)$. It 
coincides with that induced by $Z_N$
acting on the $p$-adic moduli problem (cf. [H3,§1] and [Go] for details).

Moreover

\noindent (1.1.4)  The {\it differential operator} $d = q\frac{d}{dq}$,
operating on the
$q$-expansion

\noindent sends classical forms to $p$-adic modular forms.

\medskip

  \noindent {\bf (1.2)} {\it Review of Hecke algebras.} We retain the
notation
from the previous section.

\proclaim{(1.2.1)  Definition}  The {\it Hecke algebra}
$${\bold h}  = {\bold h}(Np^{\infty},\Cal O) \subset End_{\Cal
O}(\overline{S}(N,\Cal O)) $$
is the $\Cal O$-subalgebra generated by the Hecke operators $T(\lambda)$
for all
primes $\lambda$ and by $T(\lambda,\lambda)$
for $\lambda$ relatively prime to $N$.   Similarly,
$${\bold h}^{\ell}(Np^{\alpha},\Cal O) \subset
End_{\Cal O}(\bigoplus_{k \leq \ell}S_{k}(Np^{\alpha},\Cal K))\cap
\Cal{O}[[q]]), $$
the {\it Hecke algebra of weight } $\ell$ {\it and level} $Np^{\alpha}$,
is the $\Cal O$-subalgebra generated by the $T(\lambda)$ for all
primes $\lambda$ and by the $T(\lambda,\lambda)$ for $\lambda$ relatively
prime to $N$.
\endproclaim

There are canonical surjections
of ${\bold h}(Np^{\infty},\Cal O)$ onto
${\bold h}^{\ell}(Np^{\alpha},\Cal O)$ for all $\ell, \alpha$,
and (cf. [H2, pp 243 ff.])
$$
{\bold h}(Np^{\infty},\Cal O) = \varprojlim {\bold
h}^{\ell}(Np^{\alpha},\Cal O). \leqno (1.2.2)
$$

Let ${\bold h}_{\ell}(Np^{\alpha},\Cal O)$ be the $\Cal{O}$-algebra
generated by the Hecke operators acting on $S_{\ell}(Np^{\alpha},\Cal{O})$
(note the difference between ${\bold h}_{\ell}(Np^{\alpha},\Cal O)$
and ${\bold h}^{\ell}(Np^{\alpha},\Cal O)$).  Let $\bold{e} \in \bold{h}$
denote Hida's ordinary idempotent [H2, (1.17 b)], and
let ${\bold h}^o = \bold{e}{\bold h}$ be the
universal ordinary $p$-adic Hecke algebra of level $N$.
We set $S_{\ell}^{o}(Np^{\alpha},\Cal O) =
\bold{e}S_{\ell}(Np^{\alpha},\Cal O)$,
for $\alpha \leq \infty$, and let ${\bold h}^o_{\ell}(Np^{\alpha},\Cal O)
= \bold{e}{\bold h}_{\ell}(Np^{\alpha},\Cal O)$.

Let $\gamma = 1 + Np \in Z_N$; then $X = \gamma - 1$ is
a topologically nilpotent element in
${\boldsymbol{\Lambda}}_N = \varprojlim_{\alpha}
\Cal{O}[(\Bbb{Z}/Np^{\alpha}\Bbb{Z})^{\times}].$
Similarly, one lets ${\boldsymbol{\Lambda}}=\Cal{O}[[\Bbb{Z}_p^\times]]$.
  Let ${\Lambda}$ be the Iwasawa algebra
$\Cal{O}[[X]]
\subset {\boldsymbol{\Lambda}}_N $.  The map
$$
\lambda \mapsto <\lambda>
$$
with the right hand side defined as in (1.1.3), makes ${\bold h}$
and thus ${\bold h}^o$  into continuous $\Lambda$-modules, and
one of Hida's main theorems is that

\proclaim{(1.2.3) Theorem} The Hecke algebra
${\bold h}^o$ is {\it free of finite
rank} over $\Lambda_{\Cal O}$.  Moreover, let $\ell \geq 2$, $\ell\equiv
1\,\pmod{p-1}$;
   let $P_{\ell}$
denote the element $(1+X) - (1+Np)^{\ell} \in \Lambda_{\Cal O}$.
Then the projection defines a canonical isomorphism
$$
{\bold h}^o/P_{\ell}{\bold h}^o \isoarrow  \bold{h}_{\ell}^{o}(Np,\Cal O).
$$
\endproclaim

\noindent {\bf (1.2.4) Remark.}  Note that $\Lambda_{\Cal
O}/P_{\ell}\Lambda_{\Cal O}$
is canonically isomorphic to ${\Cal O}$ for every integer $\ell \geq 1$ and
every
${\Cal O}$.

Let $\overline{M}^o = \overline{M}^o(N,\Cal O) = \bold{e}\cdot \overline{M}
\subset
\overline{M}$, $\overline{S}^o = \overline{S}^o(N,\Cal O) = \bold{e}\cdot
\overline{S}
  \subset \overline{S}$.  We define an $\Cal O$-bilinear pairing
$$
<\bullet,\bullet>: \overline{S}^o \times {\bold h}^o \rightarrow \Cal{O};
\quad <g,T> = a(1,g|T) \leqno (1.2.5)
$$
where, as usual, we write $g \mapsto g|T$ for the action of the Hecke
operator $T$.
Then (1.2.5) is a perfect pairing, with respect to which the action of
${\bold h}^o$ is
(tautologically) symmetric.  For each $\ell \geq 2$, we obtain by
restriction a perfect
pairing
$$
<\bullet,\bullet>: S_{\ell}^{o}(Np,\Cal{O}) \times
\bold{h}_{\ell}^{o}(Np,\Cal{O}) \rightarrow \Cal{O} \leqno (1.2.6)
$$
defined by the same formula as (1.2.5).  The pairing (1.2.6)
identifies $S_{\ell}^{o}(Np,\Cal{O})$
as the $P_{\ell}$-torsion submodule of $\overline{S}^o$.

\bigskip

\noindent {\bf (1.3)} {\it Congruence modules and $\Lambda$-adic forms.}

Henceforward we assume $\Cal {O}$ to be the ring of integers of a finite
extension $K$ of $\Bbb{Q}_p$.  We let $\Cal L$ denote the fraction field of
$\Lambda$,
$\Cal L^\prime$ a finite extension of $\Cal L$, and let $\bold{I}$ be the
integral
closure of $\Lambda$ in $\Cal L^\prime$.  We denote by
$\Cal{X}(\bold{I})$ the set of
prime ideals of $\bold{I}$ of height $1$, and let
$$
\Cal{X}_k(\bold{I}) = \{P \in \Cal{X}(\bold{I}) | P \cap \Lambda_{\Cal{O}}
= P_k \} .$$
For $P \in \Cal{X}(\bold{I})$ let $\Cal{O}_P = \bold{I}/P$; we let $k(P) =
k$ if
$P \in \Cal{X}_k(\bold{I})$.

  Let $\tau = \tau_{\bold{f}}: {\bold h}^o \rightarrow \bold{I}$ be a
homomorphism of $\Lambda_{\Cal{O}}$-algebras.  The subscript $\bold{f}$
refers
implicitly to a family of $p$-adic modular forms, or to a $\Lambda$-adic
modular form, in the sense of Wiles, cf. [H4].

\proclaim{(1.3.1)  Definition} We say
$\tau_{\bold{f}}$, or ${\bold{f}}$, is $N$-{\it primitive} if, for some
integer $k > 1$
(equivalently, for all $k > 1$) and for some $P \in \Cal{X}_k(\bold{I})$,
the homomorphism
$$
\tau_{\bold{f}} (mod P_{k}): {\bold h}^o/P_{k}{\bold h}^o \rightarrow
{\Cal O}_P = \bold{I}/P
$$
(cf. Remark (1.2.4)) is the homomorphism associated to a modular form of
weight
$k$ of level dividing $Np$ which is primitive at all primes dividing $N$.
\endproclaim

From now on, we consider such an $N$-primitive form ${\bold{f}}$.

\proclaim{(1.3.2) Proposition {\rom [H2, Cor. 3.7]}}   If ${\bold{f}}$ is
primitive, then
the homomorphism
$\tau_{\bold{f}}\otimes 1_{\Cal L^{\prime}}:
{\bold h}^o \otimes_{\Lambda}\Cal{L}^{\prime} \rightarrow
\Cal{L}^{\prime}$
is split over $\Cal{L}^{\prime}$.
\endproclaim

Thus there are an idempotent
$1_{\bold{f}} \in {\bold h}^o\otimes_{\Lambda}\Cal{L}^{\prime}$
and an isomorphism
$1_{\bold{f}}\cdot{\bold h}^o\otimes_{\Lambda}\Cal{L}^{\prime}
\isoarrow \Cal{L}^{\prime}$
such that the homomorphism $\tau_{\bold{f}}\otimes 1_{\Cal L^{\prime}}$ is
given by
multiplication by $1_{\bold{f}}$.  We write
$$
{\bold h}^o \otimes_{\Lambda}\Cal{L}^{\prime} \isoarrow
\Cal{L}^{\prime} \oplus \Cal{B}; \leqno(1.3.3)
$$
then $\tau_{\bold{f}}$ corresponds to projection on the first factor.

Let $\bold{h}^o_{\bold I} = \bold{h}^o
\otimes_{\Lambda}\bold{I}$,
$h(\Cal{B}) = im(\hoI) \subset \Cal{B}$ with respect to the second
projection
in (1.3.3).  It follows from Theorem (1.2.3) that (1.3.3) induces an
injection
$\hoI  \hookrightarrow \bold{I} \oplus h(\Cal{B})$ of $\bold{I}$-modules.
The {\it congruence module}
$$
\Cal{C} = (\bold{I} \oplus h(\Cal{B}))/\hoI.  \leqno (1.3.4)
$$
is an $\bold{I}$-torsion module.  On the other hand, if we take the natural
embedding
$$
i:  \bold{I} \hookrightarrow \bold{I} \oplus h(\Cal{B}); \quad \lambda
\mapsto (\lambda,0)
$$
then $\Hf = \bold{I} \cap i^{-1}(\hoI) \subset {\bold I}$ is an
ideal in $\bold{I}$.  Then $i$ induces an isomorphism
$\Cal{C} \isoarrow \bold{I}/\Hf$ of $\bold{I}$-modules.  Thus
$$
Denominator(1_{\bold{f}}) = \{a \in \bold{I} | a 1_{\bold{f}} \in \hoI \}
= \Hf.  \leqno (1.3.5)
$$
Once and for all we fix an element $H \in \Hf$ and define
$$
T_{\bold{f}} = T_{\bold{f},H} = H\cdot 1_{\bold{f}} \in \hoI.
$$

\noindent{\bf Remark:}

By assumption, $\bold{f}_{P}^{prim}$ is $N$-primitive.
Let $\omega: (\ZZ/p\ZZ)^{\times} \rightarrow \Cal{O}$
be the Teichmuller character.  Write
$$\bold{h}^o = \Pi_{a \in \ZZ/(p-1)\ZZ} \bold{h}^o(\omega^a),$$
where $\bold{h}^o(\omega^a) \subset \bold{h}^o$ is the subalgebra on
which $(\ZZ/p\ZZ)^{\times} \subset Z_N$ acts via the $a$th power $\omega^a$
of $\omega$; let $\hoI(\omega^a) = \bold{h}^o(\omega^a)
\otimes_{\Lambda} \bold{I}$.
Since $\bold{I}$ is assumed to be an integral domain, it follows
that $\tau_{\bold{f}}: \hoI \rightarrow \bold{I}$ factors through
the natural projection on $\hoI(\omega^a)$, for exactly one $a =
a(\bold{f})$, say.
For any $k \geq 2$ and any $a \in \ZZ/(p-1)\ZZ$ there is an isomorphism
$$
\bold{h}^o(\omega^a)/P_k \bold{h}^o(\omega^a) \isoarrow
\bold{h}^o_k(Np,\omega^{a-k},\Cal{O}),$$
where $\bold{h}^o_k(Np,\omega^{a-k},\Cal{O}) \subset
\bold{h}^o_k(Np,\Cal{O})$
is the subalgebra on which the quotient
$\frac{\Gamma_1(N) \cap \Gamma_0(p)}{\Gamma_1(Np)} \simeq
(\ZZ/p\ZZ)^{\times}$
acts via $\omega^{a-k}$ [H2, p.249].
If $P \in \Cal{X}_k(\bold{I})$,
  $\bold{f}_{P}^{prim}$ is thus unramified at $p$
if and only if $a \equiv k (\text{mod } p-1)$, unless $k = 2$, in which
case
the $p$-component of the automorphic representation attached to $\bold{f}$
may
be special; cf. [H1, II, p. 37].

\medskip

\noindent {\bf (1.4)} {\it Arithmetic p-adic measures.}   Let
$Z$ be a $p$-adic manifold
of the form $\Bbb{Z}_p^r \times$ {\it (finite group)}.  Let $\Cal O$
be as in (1.1), and let $\Cal{C}(Z,\Cal{O})$
denote the space of continuous $\Cal O$-valued functions on $Z$.
For any subring $R \subset \Cal{O}$ let $L\Cal{C}(Z,R)$ denote the
space of locally constant $R$-valued functions on $Z$.

Let
$$Mes(Z,\overline{M}) = Hom_{\Cal{O}}(\Cal{C}(Z,\Cal{O}), \overline{M})$$
be the set of $p$-adic measures
on $Z$ with values in $\overline{M}$.
The measure $\mu \in Mes(Z,\overline{M})$ is called {\it ordinary} (resp
{\it cuspidal})
if it takes values in $\overline{M}^o$ (resp.in $\overline{S}$).
As usual, we write
$\int_Z \phi d\mu$ in place of $\mu(\phi)$.  We assume $Z$
given with a continuous $Z_N$-action, denoted
$(z,x) ~~~ \mapsto ~~~z\cdot x$, for $z \in Z_N, x \in Z$.  We let
$z \mapsto z_p$ be the projection of $z \in Z_N$ on its $p$-adic part
$z_p \in \ZZ^{\times}_p$.

\proclaim{(1.4.1) Definition}  Let $\kappa$ be
an integer and $\xi :Z_N \rightarrow {\Cal
O}^\times$ a character of finite order.  A $p$-adic measure on $Z$
with values in $\overline{M}$ is {\it arithmetic with character }
$(\kappa,\xi)$ (cf. [H1, II, (5.1)] if

(a) for all $\phi \in L\Cal{C}(Z,\Cal{O}\cap \overline{\Bbb{Q}})$
$$
\int_Z \phi d\mu \in M_{\kappa}(Np^{\infty},\overline{\Bbb{Q}});
$$

(b)  For all $\phi$ as above
$$
(\int_Z \phi d\mu)|z = z_{p}^{\kappa} \xi(z)\int_Z \phi|z d\mu
$$
where $\phi|z(x) = \phi(z\cdot x)$ for all $z \in Z_N$.

(c)  Let $d$ denote the operator of (1.1.4).  There is
a continuous function $\nu: Z \ra \Cal{O}$ such that
$$\nu\vert z = z_p^2\cdot \nu$$
for all $z\in Z_N$ and such that, for any $\phi$ as above,
$$d^r(\int_Z \phi d\mu)=\int_Z \nu^r\phi d\mu.$$
\endproclaim

Let $\Cal{O}_{-\kappa,\xi^{-1}}$ be $\Cal{O}$, viewed as 
$\Cal{O}[[Z_N]]$-module
via the linear extension of the character $z \mapsto z_{p}^{-\kappa}
\xi^{-1}(z)$. Condition (b) can be rephrased:

{\sl (b$^{\prime}$) The map
$$\int_Z \bullet~ d\mu:\Cal{C}(Z,\Cal{O}) \rightarrow \overline S^0\otimes
\Cal{O}_{-\kappa,\xi^{-1}}$$ is ${\boldsymbol{\Lambda}}_N$-linear, 
where ${\boldsymbol{\Lambda}}_N$ acts on
$\Cal{C}(Z,\Cal{O})$
by $\Cal{O}$-linear extension of the action $\phi \mapsto \phi|z$.}

\noindent{\bf Comments:} Let us give a useful reformulation of this 
definition when
$Z =\ZZ_p^\times$. The action of $Z_N$ on $\overline{M}$ given by 
$z\mapsto <z>$
endows $\bold{h}^o$ with a natural ${\boldsymbol{\Lambda}}_N$-algebra 
structure. Let
${\boldsymbol{\Lambda}}=\Cal{O}[[\bold{Z}_p^\times]]$ and
$$\psi:Z_N\rightarrow {\boldsymbol{\Lambda}}^\times, \quad z\mapsto
z_p^\kappa\xi(z)[z_p]^2$$ For any ${\boldsymbol{\Lambda}}_N$-algebra $A$, let
$$A(\psi)=A\otimes_{{\boldsymbol{\Lambda}}_N}{\boldsymbol{\Lambda}}$$
where $ {\boldsymbol{\Lambda}}_N\rightarrow {\boldsymbol{\Lambda}}$ 
is the natural extension of
$\psi$. We view $A(\psi)$ as a ${\boldsymbol{\Lambda}}_N$-algebra.
Then, the group of ordinary cuspidal measures on $\bold{Z}_p^\times$ 
with character
$(\kappa,\xi)$ can be identified with the $\psi$-isotypic
${\boldsymbol{\Lambda}}$-submodule of $(\overline{S}^o)^\psi$ of 
$\overline{S}^o\otimes
{\boldsymbol{\Lambda}}$ defined by
$$(\overline{S}^o)^\psi= \bigcap_{z\in \bold{Z}_p^\times} 
(\overline{S}^o\otimes
{\boldsymbol{\Lambda}})^{z\otimes 1=1\otimes\psi(z)}$$
As an example, we can take the following

\noindent {\bf (1.4.2)  Example}.  Here $Z =\ZZ_p^\times$ and
$\nu$ in (c) is the tautological inclusion $Z \hookrightarrow \Cal{O}$.
The action of $Z_N$ on $Z$ is given by
$$z\cdot x=z_p^2x \quad z\in Z_N, x\in Z$$
(note the square!).  Let $R$ be the ring
of algebraic integers in $\Bbb C$, and let $g \in S_{\ell}(N,R)$ for some
$\ell$ and some
$N$ relatively prime to $p$.   We assume $g$ has
nebentypus character $\xi$.  Let $g_p$ be the twist of $g$ by the trivial
character $\pmod{p}$:
If $g=\sum_{n=1}^\infty a(n,g)q^n$,
we have $g_p=\sum_{(n,p)=1}a(n,g)q^n$.
For any function $\phi \in L\Cal{C}(Z,R)$ set
$$
\int_Z \phi d\mu_g =  \sum_{(n,p)=1} \phi(n)a(n,g) q^n.  \leqno (1.4.2.1)
$$
Extend this by continuity to $\Cal{C}(Z,\Cal{O})$ for varying $\Cal{O}$.
Hida has verified (cf. [H1, I, Prop. 8.1]) that $d\mu_g$ is arithmetic with
character $(\ell,\xi)$; its ``moments" are given by
$$
\int_Z x^r d\mu_g = d^r(g_p); \leqno (1.4.2.2)
$$
here and in what follows we write $x^r = \nu(x)^r$.

\noindent {\bf (1.4.3)  Example}.  Here $Z, R, N,$ and $g$ are as in
(1.4.2), and
we let $h \in S_m(N,R)$ for some $m$.  We let $\xi_g$ and $\xi_h$ denote
the nebentypus characters of $g$ and $h$, respectively.
The measure $h \cdot d\mu_g$ is defined by
$$
\int_Z \phi ~h \cdot d\mu_g =  h \cdot \sum_{n=1}^{\infty} \phi(n)a(n,g)
q^n.  \leqno (1.4.3.1)
$$
It follows from (1.4.2) that $h \cdot d\mu_g$ is arithmetic with character
$(\ell + m,\xi_g\cdot\xi_h)$; its
``moments" are then given by
$$
\int_Z x^r ~h \cdot d\mu_g = h \cdot d^r(g_p). \leqno (1.4.3.2)
$$

\noindent{\bf (1.4.4)} {\it Contractions of arithmetic measures by Hida
families.}
Unless otherwise specified, we assume

{(\bf G)} \quad  $Z$ is a $p$-adic group containing $\ZZ_p^\times$ as
an open subgroup
of finite index, and with an action of $(\ZZ/N\ZZ)^{\times}$.
\smallskip

\noindent Think for instance of $Z=\ZZ_p^\times$
  with trivial action of $(\ZZ/N\ZZ)^{\times}$, as we will assume
in the next chapter.

The identification of the completed group algebra ${\Cal O}[[Z]]$ with the
space of
continuous $\Cal{O}$-valued distributions on $X$, as for example in [H4],
yields an isomorphism
$${\Cal O}[[Z]] \isoarrow  Hom_{\Cal{O}}(\Cal{C}(Z,\Cal{O}),\Cal{O}).$$
By extension of scalars, we may thus identify
$$
\matrix \{\text{ ordinary cuspidal measures}  \\ \text{ on } Z \text{ with
values in } \overline{M}\} \endmatrix   \cong
Hom_{\Cal{O}}(Hom_{\Cal{O}}({\Cal O}[[Z]],\Cal{O}),\overline{S}^o).  \leqno
(1.4.4.1)$$
Let $Hom_{\Cal{O}}(Hom_{\Cal{O}}({\Cal
O}[[Z]],\Cal{O}),\overline{S}^o)_{\kappa,\xi}$
be the $\Cal{O}$-submodule of the right-hand side of (1.4.4.1)
corresponding
to arithmetic measures with character $(\kappa,\xi)$.

Let $\I$ and $\overline{S}^o_{\hat{\I}}$ be as in \S 1.3 and
let $(\kappa,\xi)$ be as in Definition 1.4.1.   Let
$Mes(Z,\overline{S}^o)_{\kappa,\xi}=(\overline{S}^o)^\psi$ be the set 
of arithmetic measures on $Z$
with character $(\kappa,\xi)$.  Suppose
\smallskip

The pairing (1.2.5) induces by extending the scalars to 
${\boldsymbol{\Lambda}}$, a pairing
$$<\cdot,\cdot>:~~~(\overline{S}^o\otimes
{\boldsymbol{\Lambda}})\hat{\otimes}_{{\boldsymbol{\Lambda}}}(\bold{h}^o\otimes
{\boldsymbol{\Lambda}})\rightarrow {\boldsymbol{\Lambda}}$$
hence

$$
Mes(Z,\overline{S}^o)_{\kappa,\xi}
\otimes_{{\boldsymbol{\Lambda}}_N(\psi)}\bold{h}^o(\psi)\rightarrow
{\boldsymbol{\Lambda}}_N(\psi)$$ We base change it to 
$\bold{h}^o(\psi)$ and obtain

$$
Mes(Z,\overline{S}^o)_{\kappa,\xi}
\otimes_{{\boldsymbol{\Lambda}}_N(\psi)}\bold{h}^o(\psi)\otimes_{{\boldsymbol{\Lambda}}_N(\psi)}\bold{h}^o(\psi)
\rightarrow \bold{h}^o(\psi)\leqno (1.4.4.3)$$
We now twist the Hida family $\tau_{\bold{f}}:\bold{h}^o\rightarrow 
\bold{I}$ by $\psi$. We
thus obtain an $\bold{h}^o(\psi)$-algebra $\bold{I}(\psi)$. We use 
this algebra to base change
$(1.4.4.3)$.
We get

$$
Mes(Z,\overline{S}^o)_{\kappa,\xi}
\otimes_{{\boldsymbol{\Lambda}}_N(\psi)}\bold{h}^o(\psi)\otimes_{{\boldsymbol{\Lambda}}_N(\psi)}\bold{I}(\psi)
\rightarrow \bold{I}(\psi)\leqno (1.4.4.4)$$

\proclaim{(1.4.5)  Definition}  Let
$$\ell_{\bold{f}}: (\overline{S}^o)^{\psi}
\rightarrow \bold{I}(\psi)$$
be the  $\boldsymbol{\Lambda}_N(\psi)$-linear map given by
$$\mu \mapsto <\mu,T_{\bold{f}}\otimes 1>.$$
We call it the contraction against the Hida family ${\bold{f}}$.
\endproclaim

%In our setting, we take $h\cdot d\mu_g\in Mes(Z,\overline{S}^o)_{\kappa,\xi}$
%and $T_\bold{f}\otimes 1\in
%(\bold{h}^o\otimes_{{\boldsymbol{\Lambda}}_N}\bold{I})(\psi)$. We obtain
%$$\theta_H(\bold{f},g,h)=(h\cdot d\mu_g)\otimes (T_\bold{f}\otimes 1)\in
%(\overline{S}^o)^\psi\hat{\otimes}_{{\boldsymbol{\Lambda}}_N(\psi)}\bold{h}^o(\psi)
%\hat{\otimes}_{{\boldsymbol{\Lambda}}_N(\psi)}\bold{I}(\psi)$$

Applying $(1.4.5)$ to $\mu = h\cdot d\mu_g ~\in 
Mes(Z,\overline{S}^o)_{\kappa,\xi}
= (\overline{S}^o)^{\psi}$
%$\theta_H(\bold{f},g,h)$, we obtain
the desired element
$$D_H(\bold{f},g,h) \in \bold{I}(\psi).$$

\noindent{\bf (1.5)}{\it Evaluation of $D_H(\bold{f},g,h)$ at certain
arithmetic points.}
\smallskip

\noindent{\bf (1.5.1)}  {\it Notation.}
\smallskip

For any form $h \in S_{k}(Np,\Qbar)$,
let $h^{\rho} = \sum_{n=1}^{\infty} \overline{a}(n,h) q^n$, where
$z \mapsto \overline{z}$ denotes complex conjugation; $h^{\rho}$ is also an
element
of $S_{k}(Np,\Qbar)$.  We let
$$\tilde{h} = h^{\rho} |_k \pmatrix 0&1\\ -Np&0 \endpmatrix.  \leqno
(1.5.1.1)$$

Let us denote by
$<\bullet,\bullet>_{p^m,k}$ the Petersson inner product for
$S_{k}(Np^m,\Cal O)$,
normalized to be linear in the first variable and anti-linear in the
second.
The formula for $<\bullet,\bullet>_{p^m,k}$ is given as usual by:
$$<f_1,f_2>_{p^m,k} =
\int_{\Gamma_0(p^m)\backslash \frak{H}} f_1(z) \overline{f}_2(z)
y^{k-2}dxdy
\leqno (1.5.1.2) $$
  whenever $f_1$ and
$f_2$
are modular forms of weight $k$ for $\Gamma_0(p^m)$ with same Nebentypus,
and one of the two is a cusp form.
\smallskip

\noindent{\bf (1.5.2)}  {\it Evaluation.}
In what follows, we let $Z = \ZZ_p^\times$.
We denote the set of height 1 prime ideals of $\bold{I}$ by
$\Cal{X}(\bold{I})$.
  Any element of $\bold{I}(\psi)$ defines
a function on $\Cal{X}(\bold{I})$.

Let $\bold{f}$ (or $\tau_\bold{f}$) be as in 1.3,
with $\Cal{O}^\prime = \Cal{O}$, and let
$g \in S_{\ell}(N,R)$ and $h \in S_{m}(N,R)$ be as in (1.4.3), where $R$ is
the ring of integers in some number field, which we assume
contained in $\Cal{O}$.
%We may regard $\ell_{\bold{f}}$ as an element of
%$Hom_{\Cal{O}}(\overline{S}^o_{\hat{\bold{I}}},\Cal{O})$; recall
%that $\ell_{\bold{f}}$ depends on the choice of an element $H \in \Hf$.
%On the other hand,
The ordinary projection
$\bold{e}(h \cdot d\mu_g)$ of the measure $h \cdot d\mu_g$ is
naturally an ordinary cuspidal measure on $Z$ of character
$(\kappa,\xi)$ for $\xi = \xi_g\cdot \xi_h$ and $\kappa = \ell + m$
%, hence
%an element of the right-hand side of (1.4.4.1).

We shall compute special values at arithmetic points of $D_H(\bold{f},g,h)$.
For $P \in \Cal{X}_k(\bold{I})$, for some $k \geq 2$, let $\bold{f}_{P}$ be the
$\bold{e}$-eigenform associated to $\bold{f}_{P}^{prim}$, with
$q$-expansion $\sum_{n=1}^{\infty} a(n,\bold{f}_{P}) q^n$, where
$$
a(np^r,\bold{f}_P) = \alpha(\bold{f}_P^{prim})^r \cdot
a(n,\bold{f}_{P}^{prim}) \text{ if } (n,p) = 1,
$$
where $\alpha(\bold{f}_{P}^{prim})$ is the $p$-adic unit root of the Hecke
polynomial
of $\bold{f}_{P}^{prim}$ at $p$.  If the nebentypus of $\bold{f}_P^{prim}$
is non-trivial
then $\bold{f}_P^{prim} = \bold{f}_P$.
In particular, $\bold{f}_{P}$ is of level exactly $Np$.

Let
$$\Cal{X}^{adm} = \{P \in \Cal{X}(\bold{I}) ~|~ \exists ~~k = ~k(P) 
~\geq ~2, ~k
\equiv 1\pmod{p-1}, ~P \in \Cal{X}_k(I)\}.$$
Then the set $\Cal{X}^{adm}$ is Zariski dense in $\Cal{X}(\bold{I})$. 
Therefore,
the element $D_H \in \bold{I}(\psi)$ is determined by its values at
points in $\Cal{X}^{adm}$.
  For any $P \in \Cal{X}^{adm}$, let $H(P) \in \Cal{O}_P$
denote the reduction of $H$ modulo $P$.
Let $P \in \Cal{X}^{adm}$;
let $T_{\bold{f},P}=H(P)\cdot 1_{\bold{f}_P}\in h_k^o(Np,\Cal{O})$,
where $k = k(P)$.   Let $2r = k - \ell - m$.
Observe that by definition of
$\psi:[z]\mapsto z_p^{\ell+m}\cdot \xi_g\xi_h(z)\cdot [z_p^2]$,
  the image
of $P$ under the twisting map $\bold{I}\rightarrow \bold{I}(\psi)$ is
above $P_{\frac{k-\ell-m}{2}} =
P_{r} \in \Lambda_{\Cal{O}}$.  Hence, by definition (1.4.5), we have

$$\aligned
D_H(\bold{f},g,h)(P) &= \ell_{\bold{f}_P}(e(h\cdot\int_Z x^r d\mu_g))\\
	&= <e(h\cdot \int_Z x^r d\mu_g)~,~T_{\bold{f},P}> \\
&= <e(h\cdot d^r g_p)~,~T_{\bold{f},P}>\\
&= \ell_{\bold{f}_P}(e(h\cdot d^r g_p))
\endaligned
\tag 1.5.2.1
$$
by compatibility of the pairings (1.2.5) and (1.2.6).

The {\it Maass operators} $\delta^r_{\ell}$,
$r = 1, 2, \dots$, defined by
Maass and Shimura, are the differential operators on the
upper half plane given by the formula
$$
\delta_{\ell} = \frac{1}{2\pi i}(\frac{\ell}{2iy} + \frac{d}{dz}); \quad
\delta^r_{\ell} = \delta_{\ell + 2r -2} \circ \cdots \circ \delta_{\ell +
2} \circ
\delta_{\ell}.$$
For any congruence subgroup
$\Gamma$, $\delta^r_{\ell}$ takes holomorphic cusp forms of weight $\ell$
for $\Gamma$ to $C^{\infty}$ functions on the upper half plane,
rapidly decreasing at infinity and ``nearly
holomorphic" in Shimura's sense [S], which transform
under $\Gamma$ like modular forms of weight $\ell + 2r$.  We refer to such
functions
as nearly holomorphic cusp forms.  If $f_i$
are nearly holomorphic cusp forms of weights $m_i$,
$i = 1, 2$, then the product $f_1 f_2$ is a nearly holomorphic cusp form of
weight
$m_1 + m_2$.

If $G$ is a nearly holomorphic cusp form of weight $k$, then
the holomorphic projection
$\Cal{H}(G)$ is the unique holomorphic cusp form of weight $k$
that satisfies
$$
<G,f>_k = <\Cal{H}(G),f>_k$$
for all holomorphic cusp forms $f$.
Let $G = h \cdot \delta^r_{\ell} g$ and $G_p = h \cdot \delta^r_{\ell}
g_p$.  It follows from
[H1, I,p. 185; II, Lemma 6.5, (iv)] that
$$
\bold{e}(h \cdot d^rg_p) = \bold{e}(\Cal{H}(G_p)).
$$
Thus, returning to formula (1.5.2.1), we find that
$$
D_H(\bold{f},g,h)(Q) =
\ell_{\bold{f},P} \circ \bold{e}(\Cal{H}(G_p)).
$$
Finally, appealing to [H1, I,prop. 4.5, II, 7.6] , we find that
$$
D_H(\bold{f},g,h)(Q) =
H(P) \cdot \alpha(\bold{f}_{P}^{prim})^{-1}\cdot p^{k-1}\cdot
\frac{<G_p,\tilde{\bold{f}}_P(pz)>_{p^2,k}}{<\bold{f}_P,\tilde{\bold{f}}_P>_
{p^2,k}}.
\leqno (1.5.2.2)
$$
\bigskip

\centerline{\bf 2. Triple product L-functions}
\medskip

\noindent{\bf (2.1)} {\it A formula for the central critical value.}
We retain the notation of the previous section.
Let $f \in S_k(N,R), g \in S_{\ell}(N,R), h \in S_m(N,R)$ be three modular
forms of
level $N$, with $k \geq \ell \geq m$.  We write their standard Hecke
$L$-functions
as follows:
$$
L(s,?) = \prod_{(q,N) = 1} L_p(s,?) \times \prod_{q|N}L_q(s,?), \quad ? =
f, g, h$$
where, for $(q,N) = 1$ the local Euler factors are of the form
$$\aligned
L_q(s,f) &= [(1 - \alpha_{1,q}q^{-s})(1 - \alpha_{2,q}q^{-s})]^{-1},\\
L_q(s,g) &= [(1 - \beta_{1,q}q^{-s})(1 - \beta_{2,q}q^{-s})]^{-1},\\
L_q(s,h) &= [(1 - \gamma_{1,q}q^{-s})(1 - \gamma_{2,q}q^{-s})]^{-1},
\endaligned \leqno (2.1.1)$$

Here our $L$-functions are normalized so that
$|\alpha_{i,q}| = q^{\frac{k-1}{2}}, |\beta_{i,q}| = q^{\frac{\ell-1}{2}},
|\gamma_{i,q}| = q^{\frac{m-1}{2}}, i = 1, 2$, for any
archimedean absolute value.  The triple product
$L$-function is the convolution of these three:
$$
L(s,f,g,h) = \prod_{(q,N) = 1}[\prod_{i,i^{\prime},i^{\prime \prime} = 1,
2} (1 - \alpha_{i,q}\beta_{i^{\prime},q}\gamma_{i^{\prime
\prime},q}q^{-s})]^{-1} \times \prod_{q|N} L_q(s,f,g,h)  \leqno (2.1.2)
$$
where the factors $L_q(s,f,g,h)$ for $q|N$ are the local Artin $L$-factors
of the corresponding
Weil-Deligne group representations, defined by reference to
the local Langlands correspondence for $GL(2)$.

In what follows we restrict attention to the case $N = 1$, i.e.,
we assume our forms are all of level $1$.  This implies in particular that
the weights $k, \ell, m$ are all {\it even}.  We assume that
$f$, $g$, and $h$ correspond to cuspidal automorphic representations
$\pi(f), \pi(g)$, and $\pi(h)$, respectively, of $GL(2)_{\QQ}$, with
trivial central characters $\xi(f), \xi(g), \xi(h)$, respectively.
denote the respective central characters.

The analytic continuation of
the triple product $L$-function has been proved by the method of
Langlands-Shahidi
[Sha] and by a variant of the Rankin method, due to Garrett [G1,G2] and
generalized by Piatetski-Shapiro and Rallis [PSR].  It is known
to satisfy a functional equation of the usual type, relating the values
at $s$ and $w + 1 - s$, where $w = k + \ell + m - 3$.

We assume henceforward that
$$ k \geq \ell + m. \leqno (2.1.3)$$
Under hypothesis (2.1.3), a formula is obtained in [HK1] -- the Main
Identity 9.2 --
for the central critical value $L(\frac{w+1}{2},f,g,h)$ of the triple
product
$L$-function.  The value is expressed as an integral of theta lifts of
$f, g, h$ to the orthogonal group attached to the split quaternion algebra
$M(2)_{\QQ}$ over $\QQ$; i.e. to the split form of $O(4)$.
The exact formula depends on several auxiliary choices.  Let $H$ denote the
algebraic group
$(GL(2) \times GL(2))/d(\Bbb{G}_m)$, where $d$ is the diagonal embedding.
Then $H$ is naturally isomorphic
to the identity component of the group of orthogonal similitudes of
the split quaternion algebra.
We let $\Cal{S}$ be the space of
Schwartz-Bruhat functions on $M(2)(\bold{A})$.
To any $\phi \in \Cal{S}$ that satisfies appropriate
finiteness properties ($K$-finite for a maximal compact subgroup
$K$ of $GL(2,\RR) \times O(2,2)$ with respect to the Weil representation)
and any automorphic form $F$ on $GL(2,\QQ)\backslash GL(2,\bold{A})$, the
theta
correspondence associates
an automorphic form $\theta_{\phi}(F)$ on $H(\QQ)\backslash H(\bold{A})$
(see [HK2], (5.1.12) for the precise formula, which
also depends on the choice of a measure, specified in [HK1]).

  Let $r = \frac{k - \ell - m}{2}$, which by (2.1.3) is a
positive integer.  Let $f^{\iota}$ be the normalized newform whose Hecke
eigenvalues are the complex conjugates of those of $f$.  In fact
$f^{\iota} = f$, since $N = 1$, but we leave the notation $f^{\iota}$
with a view to future generalizations.  Then $\overline{f}^{\iota}$ is
an antiholomorphic form with the same Hecke eigenvalues as $f$; in other
words, $\overline{f}^{\iota}$ lifts to an element $f^!$ of $\pi(f)$.
Similarly, let $g^!(i)$ and $h^!$ be liftings of
$ \delta^i_{\ell}(g), ~~0 \leq i \leq r$ and $h$,
respectively, to automorphic forms on
$GL(2,\QQ) \RR^{\times}_{+}\backslash GL(2,\bold{A})$; i.e., to elements of
$\pi(g)$ and $\pi(h)$, respectively.  Here we are using the fact
that the Maass operators correspond to elements of $Lie(GL(2))$, cf.
[HK1, Lemma 12.5] and the references cited there.  It follows
from [HK1, Theorem 7.2] that, for appropriate choices of
$\phi^i \in \Cal{S}$, $i = 1, 2, 3$, we can arrange that
$$
\theta_{\phi^1}(f^!) = f^!;~ \theta_{\phi^2}(g^!(0)) = g^!(r);~
\theta_{\phi^3}(h^!) = h^! \leqno (2.1.4)$$

We abbreviate $\Phi = (\phi^1,\phi^2,\phi^3), F = (f^!, g^!(r),h^!)$.
Let $d\mu$ be the $GL(2,\bold{A})$-invariant Haar measure on
$\bold{A}^{\times}\cdot GL(2,\QQ) \backslash GL(2,\bold{A})$ with total
volume $1$.  Then we have the following formula:

\proclaim{MAIN IDENTITY}  ([HK1, 9.2]):
$$
Z_{\infty}(F,\Phi) \cdot\cdot L(\frac{w+1}{2},f,g,h) =
2\zeta(2)^2 \cdot I(f^!, g^!(r), h^!)^2, \leqno (2.1.5)
$$
where
$$
I(f^!, g^!(r), h^!) =
\int_{\bold{A}^{\times}\cdot GL(2,\QQ) \backslash GL(2,\bold{A})} f^! \cdot
g^!(r) \cdot h^! d\mu
\leqno (2.1.6)
$$
\endproclaim
Here $\zeta(2)$ is the value at $s = 2$ of the Riemann zeta function, and
$Z_{\infty}(\bullet,\bullet)$ is the value
at $s = 0$ of the normalized local zeta
integrals, defined by Garrett and Piatetski-Shapiro-Rallis.  The nature
of $Z_{\infty}(\bullet,\bullet)$ will be discussed in the next section;
here we
merely remark that the notation of [HK1] has been slightly simplified in
the present account.

\medskip

\noindent{\bf (2.2)} {\it $p$-adic interpolation of certain central
critical values.}
\smallskip

The Main Identity (2.1.6) can be rewritten
$$
(\frac{<h \cdot \delta^r_{\ell}(g),f>_k}{<f,f>_k})^2 =
\frac{Z_{\infty}(F,\Phi)}{2\zeta(2)^2}
\frac{L(\frac{w+1}{2},f,g,h)}{(<f,f>_k)^2}. \leqno (2.2.1)$$

The left hand side of (2.2.1)
has almost the same form as the
square of a special value (1.5.2.2)
of the $p$-adic measure constructed in 1.5.  The only
modification necessary is to replace $f$ by the value $\bfP$ at a prime $P$
of
an ordinary Hida family, and to incorporate the twist $\bfP \mapsto \bftP$.

Write $G = h \cdot \delta^r_{\ell}(g)$ and $G_p = h \cdot
\delta^r_{\ell}(g_p)$, as in \S 1.5.2.  Let $\bold{f}$,
$\bfP$, and $\bfpP$ be as in section 1.3. In what follows, we let $f =
\bfpP$.
Recall that we have fixed the auxiliary level $N$ to be $1$. Let
$\alpha_1 = \alpha_1(\bfP)$ be the $p$-adic unit root of the Hecke
polynomial
of $\bfpP$ at $p$, and let $\alpha_2 = \alpha_2(\bold{f}_P)$ denote its
other
root.  Recall that $<\bullet,\bullet>_{p^m,k}$ ($m\geq 0$) has been defined
in (1.3.11).

\proclaim{Proposition 2.2.2}  With notations as above, the following
formula is valid:

$$\frac{<G_p,\bftP>_{p^2,k}}{<\bfP,\bftP>_{p,k}} =
\frac{E_p(\bfP,g,h)}{p^{1-\frac{k}{2}}\alpha_1 (1 -
\frac{\alpha_2}{\alpha_1})(1 -
\frac{\alpha_2}{p\alpha_1})}
\cdot \frac{<G,\bfpP>_{1,k}}{<\bfpP,\bfpP>_{1,k}}.$$

where
$$E_p(\bfP,g,h) =p^{-k}(p^2\alpha_1^2-\alpha_2
a_p)-p^{2-\frac{k+\ell+m}{2}}\alpha_1 b_p c_p
+p^{1-\frac{k+\ell-m}{2}}b_p^2+p^{1-m}c_p^2-\alpha_2
p^{1-\frac{k+m}{2}}c_p-1
$$
\endproclaim

\demo{Proof}  The elements of this calculation
are certainly well known to specialists.  However,
we were unable to find a complete comparison of the two sides
in the literature, so we are including all details.

We extend the Petersson inner product $\langle \phi,\psi\rangle_{p^m,k}$ to
$C^\infty$ forms
of and weight $k$, level
$p^m$ with trivial Nebentypus, one of them decreasing rapidly at cusps.
Recall that
$\f_P=f-\al_2\cdot f\vert [p]$ where:
\medskip

$\quad \quad \quad \bullet \quad a(p,f)=\al_1+\al_2$,
$\al_1\al_2=p^{k-1}$ and $\al_1$ is a $p$-adic unit, and
\smallskip
$$\bullet \quad \phi\vert[m]=\phi(mz)=m^{-{\ell\over 2}}\cdot \phi\vert
\pmatrix m&0\\0&1\endpmatrix)  \qquad \qquad \qquad \qquad \qquad$$
for any $\phi$ of weight $\ell$ and any $m\geq 1$.

  From
$$\widetilde{\f}_P=\f_P\vert\pmatrix
  0&-1\\p&0\endpmatrix$$
using the equality
$\pmatrix 0&1\\ -p&0 \endpmatrix = \pmatrix 0&-1\\ 1&0 \endpmatrix \p1$,
we find
$$\widetilde{\f}_P(pz)=p^{-{k\over 2}}\cdot\left(f\vert
\pmatrix p^2&0\\0&1\endpmatrix -\al_2\cdot p^{-{k\over
2}}f\vert\pmatrix p&0\\0&1\endpmatrix\right).$$

Set $A(p^m)=\pmatrix p^m&0\\0&1\endpmatrix $.
We will repeatedly use  the following

\proclaim{Lemma 2.2.3}  The following formulas are valid:

1) Given any $C^\infty$ forms $\phi$ and $\psi$ of weight $w$, level $1$,
with $\phi$ eigen for $T_{p^m}$
of eigenvalue $\lambda_{p^m}$ :
$$\align &\text{(i)} <\phi,\psi|_w ~A(p^m)>_{p^m,w} =
p^{m(1-\frac{w}{2})}\cdot \lambda_{p^m} <\phi,\psi>_{1,w}.\\
&\text{(ii)} <\phi,\psi>_{p^m,w} = [\Gamma_0(p^m):SL(2,\ZZ)]\cdot
<\phi,\psi>_{1,w}.\\
&\text{(iii)} <\phi|_w ~A(p),\psi|_w ~A(p)>_{p,w} =
(p+1)<\phi,\psi>_{1,w} .\endalign$$

2) Similarly, if $\psi$ has level $p$  and $\phi$ level $1$ and is eigen
for $T_p$, one has:

$$\align &\text{(iv)} <\phi,\psi|_w ~A(p)>_{p^2,w} =
p^{1-\frac{w}{2}}\cdot \lambda_{p}\cdot <\phi,\psi>_{1,w}-<\phi|_w
A(p),\psi >_{p,w}.\endalign$$
\endproclaim

\demo{Proof (of Lemma 2.2.3)}  We proceed as in [PR,4.2 or H1, II, Lemma
5.3] by observing
$$SL(2,\ZZ) ~A(p^m) \Gamma_0(p^m) = SL(2,\ZZ) ~A(p^m)$$
and, if $\phi, \psi \in S_k(SL(2,\ZZ))$ and $\gamma \in GL(2,\QQ)$
has positive determinant,
$$<\phi,\psi|_w[SL(2,\ZZ)\gamma \Gamma_0(p^m)]>_{p^m,w}
= <\phi|_w [\Gamma_0(p^m) \gamma^{\iota}SL(2,\ZZ)],\psi>_{1,w}.$$
Here $\gamma \mapsto \gamma^{\iota}$ is the main involution
$\pmatrix a&b\\c&d \endpmatrix^{\iota} = \pmatrix d&-b\\-c&a \endpmatrix.$
Moreover, one checks that if
$$\Gamma_0(p^m)\pmatrix
p^m&0\\0&1\endpmatrix SL_2(\ZZ)=\coprod_i\Gamma_0(p^m)\alpha_i$$
then
$$SL_2(\ZZ)\pmatrix
p^m&0\\0&1\endpmatrix SL_2(\ZZ)=\coprod_iSL_2(\ZZ)\alpha_i$$
and since $\phi$ has level $1$, one finds
$$
p^{m(\frac{w}{2} -1)}<\phi,\psi|_k~A(p^m)>_{p^m,w} = <\phi|T_{p^m},
\psi>_{1,w}
= \lambda_p<\phi,\psi>_{1,k},$$
which proves (i).

Next, assertion (ii) is obvious, and (iii) is similar, when $\Gamma_0(p^m)$
is replaced
by $A(p^m)\Gamma_0(p^m) A(p^m)^{-1}$, which has the same index in
$SL(2,\ZZ)$.

For assertion (iv), one observes the equality of sets
$$\Gamma_0(p)A(p)=\Gamma_0(p) A(p) \Gamma_0(p^2)$$
then one uses the adjunction formula for
$$[\Gamma_0(p) A(p) \Gamma_0(p^2)]$$ together with the fact that
$\Gamma_0(p^2) A(p)^\iota \Gamma_0(p)$ and $U_p=\Gamma_0(p) A(p)^\iota
\Gamma_0(p)$ admit a same
  set of representatives.
\enddemo

\smallskip

\noindent{\bf Step 1:  Computation of} $<G_p,\bftP>_{p^2,k}$
\smallskip

We have
$$<G_p,\bftP>_{p^2,k}=p^{\frac{-k}{2}}\cdot
\langle f\vert A(p^2)-\al_2\cdot p^{-{k\over 2}}f\vert A(p),h\cdot
\delta^r_\ell(g\vert \iota_p)\rangle_{p^2,k}$$

where, if $g=\sum_{n\geq 1} b_nq^n$, one has $g\vert
\iota_p=\sum_{(n,p)=1}b_nq^n$ and
$r=k-\ell-m$.

  Let us observe that since
$g$ is a Hecke eigenform,
$$g\vert \iota_p=g\vert(1-T_p[p]+p^{\ell-1}[p^2])$$
Therefore,
$$(2.2.4)\quad h\cdot \delta^r_\ell(g\vert \iota_p)=h\cdot \delta^r_\ell
g-b_p p^{-{\ell\over
2}} h\cdot \delta^r_\ell \left( g\vert
\pmatrix p&0\\0&1\endpmatrix \right)+p^{-1}  h\cdot \delta^r_\ell \left(g\vert
\pmatrix
p^2&0\\0&1\endpmatrix\right) $$
Recall by the way that
$$\delta^r_\ell \left(g\vert_\ell \al\right)=\left(\delta^r_\ell g\right)
\vert_{\ell+2r} \al.$$
Now, by substituting (2.2.4) in the inner product
$$
\langle f\vert A(p^2)-\al_2\cdot p^{-{k\over 2}}f\vert A(p),h\cdot
\delta^r_\ell(g\vert \iota_p)\rangle_{p^2,k}$$
one obtains a sum of six terms that we
compute separately. Let $G=h\cdot \delta^r_\ell g$.
\smallskip
\noindent $\bullet$

$$T_1=\langle f\vert A(p^2),   G\rangle_{p^2,k}$$

Since $G$ has level $1$, we can apply Lemma 2.2.3 (i); one finds
$$T_1=p^{2-k}\cdot \langle f,   G\vert
T_{p^2}\rangle_{1,k}=p^{2-k}(a_p^2-p^{k-1})\cdot \langle
f,   G\rangle_{1,k}
$$
\smallskip
\noindent $\bullet$

$$T_2=-\al_2 p^{-{k\over 2}} \langle f\vert A(p),  G\rangle_{p^2,k}$$
Observe
$$\langle f\vert A(p),   G\rangle_{p^2,k}=p\cdot \langle f\vert
\pmatrix
p&0\\0&1\endpmatrix,   G\rangle_{p,k}$$
Then, by the same reasoning as above, one has
$$T_2=-\al_2 p^{-{k\over 2}}p\cdot p^{1-{k\over 2}}\cdot \langle f,
G\vert T_{p}\rangle_{1,k}$$
so,
$$T_2=-\al_2 p^{2-k}a_p\cdot \langle f,   G\rangle_{1,k}$$

\smallskip
\noindent $\bullet$
$$T_3=-b_p\cdot p^{-{\ell\over 2}}\cdot \langle f\vert A(p^2),   h\cdot
\delta^r_\ell g\vert
  A(p)
\rangle_{p^2,k}$$
  One rewrites $T_3$ as
$$-b_p\cdot p^{1-{\ell+m\over 2}}\cdot \langle \left(f\vert
\pmatrix
p&0\\0&1\endpmatrix \overline{\delta^r_\ell g}y^{k-m}\right)\vert
\Gamma_0(p)
\pmatrix p&0\\0&1\endpmatrix\Gamma_0(p^2),   h
\rangle_{p^2,m}$$
By Lemma 2.2.3 (iv), one gets
$$p^{{\frac{m}{2}-1}}\langle \phi\vert A(p),h\rangle_{p^2,m}=
\langle \phi,h\vert U_p\rangle_{p,m}$$
where $h\vert U_p=h\vert T_p-p^{{m\over 2}-1}\cdot h\vert A(p) $.
Thus, one has
$$T_3=- p^{1-{\ell+m\over 2}}b_pc_p\cdot \langle f\vert A(p)\cdot\overline{
\delta^r_\ell g}y^{k-m},h
\rangle_{p,m}+p^{-{\ell\over 2}}b_p \langle f\vert A(p)\cdot\overline{
\delta^r_\ell g}y^{k-m},h\vert A(p)
\rangle_{p,m}$$
and finally,
$$T_3=- p^{2-{k+\ell+m\over 2}}a_p b_p c_p\cdot \langle
f,G\rangle_{1,k}+p^{1-{k+\ell-m\over 2}}b_p^2\cdot \langle
f,G\rangle_{1,k}$$

\smallskip
\noindent $\bullet$

$$T_4=\al_2b_p p^{-{\ell+k\over 2}}\cdot \langle f\vert A(p),   h\cdot
\delta^r_\ell g\vert A(p)\rangle_{p^2,k}$$
We rewrite it as
$$\al_2b_pp^{-{\ell+k\over 2}}\cdot \langle\left(f\overline{\delta^r_\ell
g} y^{k-m}\right)\vert A(p),   h\rangle_{p^2,m}.$$
That is,
$$T_4=\al_2b_pc_pp^{2-{k+\ell+m\over 2}}\cdot \langle f,G\rangle_{1,k}$$

\smallskip
\noindent $\bullet$

$$T_5=p^{-1}\cdot \langle f\vert A(p^2),   h\cdot \delta^r_\ell g\vert
A(p^2)\rangle_{p^2,k}.$$
By the same calculation, we get
$$T_5=p^{1-m}(c_p^2-p^{m-1})\cdot\langle f,G\rangle_{1,k}$$

\smallskip
\noindent $\bullet$

$$T_6=-\al_2 p^{-1-{k\over 2}}\cdot  \langle f\vert A(p),   h\cdot
\delta^r_\ell g\vert  A(p^2)\rangle_{p^2,k}$$
which is equal to
$$-\al_2 p^{-1-{k\over 2}}p^{1-{m\over 
2}}\langle\left(f\overline{\delta^r_\ell g} y^{k-m}\right)\vert
  A(p),   h\rangle_{p^2,m}.$$
Hence by adjunction
$$T_6=-\al_2 p^{-1-{k\over 2}}\left[ p^{1-{m\over 2}}c_p\langle
f\overline{\delta^r_\ell g} y^{k-m},h\rangle_{p,m}
-\langle f\overline{\delta^r_\ell g} y^{k-m},h\vert
A(p^2))\rangle_{p,m}\right]
$$
Note that  $\langle f\overline{\delta^r_\ell g}
y^{k-m},h\rangle_{p,m}=p\cdot \langle f\overline{\delta^r_\ell g}
y^{k-m},h\rangle_{1,m}$
and
$$\langle f\overline{\delta^r_\ell g}\vert A(p)\cdot y^{k-m},h\vert
A(p)\rangle_{p,m}=\langle f,\left(\delta^r_\ell
  g\cdot h\right)\vert A(p) \rangle_{p,k}=p^{1-{k\over 2}}a_p\cdot\langle
f,G\rangle_{1,k}$$
Therefore,
$$T_6=-\al_2 p^{1-{k+m\over 2}}c_p\cdot \langle f,G\rangle_{1,k}+\al_2a_p
p^{-k}\cdot \langle f,G\rangle_{1,k}$$

\smallskip
The sum of the terms $T_i$ ($i=1,\ldots,6$) is the product of $ \langle
f,G\rangle_{1,k}$ by
$$ p^{2-k}(a_p^2-p^{k-1})-\al_2 p^{2-k}a_p- p^{2-{k+\ell+m\over 2}}a_p b_p
c_p+p^{1-{k+\ell-m\over
2}}b_p^2
+\al_2b_pc_pp^{2-{k+\ell+m\over 2}}+$$
$$+p^{1-m}(c_p^2-p^{m-1})-\al_2 p^{1-{k+m\over 2}}c_p+\al_2a_p p^{-k}$$
That is,
$$
\frac{<G_p,\bftP>_{p^2,k}}{<G,f>_{1,k}}=E_p(\bfP,g,h) \leqno (2.2.5)$$

\noindent{\bf Step 2:  Computation of} $<\bfP,\bftP>_{p,k}$
\smallskip

$$\align
<\bfP,\bftP>_{p,k} &= <f,f|_k A(p))>_{p,k} -
\alpha_2 \cdot p^{\frac{k}{2}}<f|_k A(p),f|_k A(p)>_{p,k}\\
&-\alpha_2 \cdot p^{\frac{k}{2}}<f,f>_{p,k}
+ \alpha_2^2 \cdot p^k<f|_k A(p),f>_{p,k}.\endalign$$

It follows from Lemma 2.2.3 that
$$
\frac{<\bfP,\bftP>_{p,k}}{<f,f>_{1,k}} =
p^{-\frac{k}{2}}(p\cdot a_p - 2(p+1)\alpha_2 + p^{1-k}\alpha_2^2\cdot
\overline{a}_p).$$
Since $f$ is of level $1$, $\overline{a}_p = a_p$ and
$\alpha_1 \cdot \alpha_2 = p^{k-1}$. Thus,

$$
\frac{<\bfP,\bftP>_{p,k}}{<f,f>_{1,k}} =
p^{1-\frac{k}{2}}\alpha_1 (1 - \frac{\alpha_2}{\alpha_1})(1 -
\frac{\alpha_2}{p\alpha_1}). \leqno (2.2.6)$$

Proposition 2.2.2 now follows immediately by combining (2.2.5) and (2.2.6).

\enddemo

Recall we put

$$ E_p(\bfP,g,h)=p^{-k}(p^2\alpha_1^2-\alpha_2
a_p)-p^{2-\frac{k+\ell+m}{2}}\alpha_1 b_p c_p
+p^{1-\frac{k+\ell-m}{2}}b_p^2+p^{1-m}c_p^2-\alpha_2
p^{1-\frac{k+m}{2}}c_p-1
$$

Let $S(P) = (1 - \frac{\alpha_2}{\alpha_1})(1 -
\frac{\alpha_2}{p\alpha_1})$.  It follows from Lemma 2.2.2
that the right-hand side of (1.5.2.2) equals
$$
H(P) \cdot \alpha_1^{-2}p^{k-2}\frac{E_p(\bfP,g,h)}{S(P)}\cdot
\frac{<G,f>_{1,k}}{<f,f>_{1,k}}.  \leqno (2.2.7)$$

Let
$$K(P)= \alpha_1^{-2}p^{k-2}\frac{E_p(\bfP,g,h)}{S(P)}$$

Combining (1.5.2.2) with (2.2.1), we then obtain our main result.

\proclaim{Theorem 2.2.8} Let $\bold{f}$ be a $p$-adic family of ordinary
cusp
forms, in the sense of 1.3, unramified outside $p$.  Let $g$ and $h$ be
cusp forms
of weights $\ell$ and $m$, respectively, of level $1$.  Let $H$ be an
annihilator
of the congruence module attached to $\bold{f}$, and let
$D_H(\bold{f},g,h)$
be the generalized $p$-adic measure constructed in section 1.4.  For any
integer $k\geq 2$, $k\equiv 1\,\pmod{p-1}$ and for $P \in 
\Cal{X}_k(\bold{I})$, the
value of
$D_H(\bold{f},g,h)$ at $P$ is related to
the central critical value $s = \frac{w+1}{2}$
of $L(s,\bfpP,g,h)$ by the following formula:
$$
(\frac{D_H(\bold{f},g,h)(P)}{H(P)\cdot K(P)})^2 =
\frac{Z_{\infty}(F,\Phi)}{2\zeta(2)^2}
\frac{L(\frac{w+1}{2},\bfpP,g,h)}{(<\bfpP,\bfpP>_k)^2}.  \leqno (2.2.9)
$$
\endproclaim

\bigskip

\noindent{\bf (2.3)} {\it Refinement of the main formula.}
\smallskip

In order to compare the results described in the main formula to accepted
conjectures
on $p$-adic $L$-functions, or to formulate reasonable conjectures regarding
the square roots of $p$-adic $L$-functions along ``anti-cyclotomic"
variables,
it would be necessary to determine the $p$-adic
behavior of the archimedean zeta integral $Z_{\infty}(F,\Phi)$
as the weight $k$ varies.  The local nature of the calculations
in [HK1] makes it clear that $Z_{\infty}(F,\Phi)$ depends only
on the weights $k, \ell, m$.  Our choices of
archimedean data are dictated by the $p$-adic construction, and a full
calculation of the archimedean integral would require summing $r = \frac{k
- \ell - m}{2}$
separate terms for given $k$.
Ikeda has recently computed archimedean triple product zeta integrals under
very general hypotheses [I].  In the case  $k \geq \ell
+ m$ his inputs are
not quite the same are ours, but his techniques may be applicable to
determine
$Z_{\infty}(F,\Phi)$ explicitly.  We note that
$Z_{\infty}(F,\Phi)$ has been determined up to rational multiples in [HK1].
Bearing in mind the slightly different normalization used in [HK1], we find
that
$$
\frac{Z_{\infty}(F,\Phi)}{\pi^{4-2k}} \in \QQ^{\times}.
$$
Since $\zeta(2) = \frac{\pi^2}{6}$,
we recover the statement of the introduction.

	We have restricted attention to forms $f$, $g$, and $h$ of level $1$.
Allowing ramification at primes different from $p$ will modify the final
formula.
We may treat bad finite places $v$ as we have treated the infinite place,
choosing local Schwartz-Bruhat functions $\phi_{i,v}$.
Then nothing will change on the right-hand side
of the Main Identity (2.1.6), but the left-hand side will include
additional zeta
integrals $Z_{v}(F,\Phi)$, reflecting these choices.  Just as in the
archimedean
case, these local zeta integrals will depend only on the local components
of
the automorphic representations associated to $f$, $g$, and $h$.  The
qualitative variation of these ramified components in a Hida family has
been determined by Hida [H3, pp. 129-133]
(the variation in a Hida family attached to Hecke
characters of an imaginary quadratic field can be seen quite explicitly).
In any case, for fixed conductor $N$, the number of distinct possible bad
non-archimedean components is finite, so the possible denominators created
by the
local integrals $Z_{v}(F,\Phi)$ remain bounded.

Our restriction to level $1$ is more serious at the prime $p$.  Requiring
that $\bold{f}_{P}$ be unramified at $p$ for all $P$ amounts to restricting
attention to a single branch of the Hida family $\bold{f}$, namely to
those $k$ congruent to $a(\bold{f})$ (mod $p-1$). In general, one
wants to allow the conductor of $\bold{f}_{P}$ to be divisible
by $p$ but not by $p^2$.  Removing the restriction that $\bold{f}_{P}$ be
unramified at $p$
only makes sense if we also allow $g$ and $h$ to have conductor divisible
by $p$.  In that case, the Main Identity will involve a zeta integral
$Z_{p}(F,\Phi)$, which might introduce additional $p$-adic zeroes or poles.
Explicit determination of such a local integral is also extremely
difficult.
The case of three special representations is considered in the article [GK]
of
Gross and Kudla; its explicit calculation is one of the most intricate
in the theory of $L$-functions.

Anyone who successfully undertakes these calculations should find it
easy, using the methods of this paper, to construct a $p$-adic measure in
three variables (allowing $f$, $g$, $h$ and the nebentypus characters
to vary, always subject to $\xi(f)\cdot \xi(g) \cdot \xi(h) = 1$), whose
moments interpolate the square roots
of normalized central critical values of triple product $L$-functions.
\bigskip
\medskip

  \centerline{\bf  REFERENCES}
\medskip
\eightpoint
\noindent [B] Bourbaki, N.: {\it Alg\`ebre Commutative}.  Paris:  Hermann,
1961, 1964,
1965, and 1983.
\smallskip

\noindent [G1] Garrett, P. B.:  Integral representations of certain
L-functions
attached to one, two, and three modular forms (manuscript, 1985).

\smallskip
\noindent [G2]  Garrett, P. B.:  Decomposition of Eisenstein series:
Rankin
triple products, {\it Annals of Math.}, {\bf 125} (1987) 209-235.
\smallskip

\noindent [Go]  Gouv\^ea, F. Q.: {\it Arithmetic of p-adic Modular Forms},
{\it Lecture Notes in Math.} {\bf 1304} (1988).
\smallskip

\noindent [GT]  Greenberg, R, and J. Tilouine: The behavior of the
symmetric
square $p$-adic $L$-function at $s = 1$, (to appear).
\smallskip

\noindent [GK]  Gross, B. and S. S. Kudla: Heights and the central critical
values of triple product $L$-functions, {\it Compositio Math.}, {\bf 81}
(1992)  143-209.

\smallskip

\noindent [HK1]  Harris, M. and S. S. Kudla: The central critical value of a
triple product $L $-function, {\it Annals of Math.}, {\bf 133} (1991)
605-672.
\smallskip

\noindent [HK2]  Harris, M. and S. S. Kudla: Arithmetic automorphic forms
for
the non-holomorphic discrete series of $GSp(2)$,
{\it Duke Math. J.}, {\bf 66} (1992) 59-121.
\smallskip

\noindent [H1]  Hida, H.:  A $p$-adic measure attached to the zeta
functions
associated with two elliptic modular forms I, {\it Invent. Math.} {\bf 79},
159-195;
II, {\it Ann. Inst. Fourier} {\bf 38} (1988)  1-83.
\smallskip

\noindent [H2]  Hida, H.:  Iwasawa modules attached to congruences of cusp
forms,
{\it Ann. Scient. E.N.S. 4-\`eme s\'erie}, {\bf 19} (1986)  231-273.
\smallskip

\noindent [H3]  Hida, H.:  $p$-adic $L$-functions for base change
lifts of $GL_2$ to $GL_3$, in L. Clozel and J.S. Milne, eds,
{\it Automorphic Forms, Shimura Varieties, and L-functions},
{\it Perspectives in Mathematics}, {\bf 11}, Vol. II, 93-142 (1990).

\noindent [H4]  Hida, H.:  {\it Elementary Theory of $L$-functions and
Eisenstein
series}, {\it London Mathematical Society Student Texts} {\bf 26},
Cambridge:  Cambridge University Press (1993).

\smallskip

\noindent [H5]  Hida, H.: On $\Lambda$-adic forms of half-integral weight
for
$SL(2)_{/\QQ}$., in S. David, ed., {\it Number Theory, Paris 1992-93},
Cambridge:  Cambridge
University Press (1995) 139-166.

\smallskip

\noindent [HTU]  Hida, H., J. Tilouine, and E. Urban:
Adjoint modular Galois representations and their Selmer groups,
{\it Proc. NAS Conference on Elliptic Curves and Modular Forms, Washington
D.C.} (to appear).

\smallskip

\noindent [I]  Ikeda, T.: On the gamma factor of the triple
$L$-function, {\it Comp. Math.},
{\bf 117} (1999); II, {\it J. Reine Angew. Math}, {\bf 499} (1998) 199-223.
\smallskip

\noindent [K]  Katz, N.: $p$-adic interpolation of real analytic
Eisenstein series, {\it Annals of Math.}, {\bf 104} (1976)  459-571.

\smallskip
\noindent [M1]  Mori, A.:  A characterization of integral elliptic modular
forms, Brandeis University thesis, (1989).
\smallskip

\noindent [M2]  Mori, A.:  A characterization of integral elliptic
automorphic
forms, {\it Ann. Sc. Norm. Sup. Pisa}, {\bf 21} (1994)  45-62.
\smallskip

\noindent [O]  Orloff, T.: Special values and mixed weight triple products,
{\it Invent. Math.}, {\bf 90} (1987)  169-180.
\smallskip

\noindent [P]  Panchishkin, A. A.:  Familles $p$-adiques de
repr\'esentations
galoisiennes et de fonctions $L$ associ\'ees aux triplets de formes
modulaires, expos\'e au S\'eminaire de th\'eorie des nombres,
Universit\'e Paris XIII, 27 janvier 1994.

\smallskip

\noindent [PR]  Perrin-Riou, B.: Fonctions $L$ $p$-adiques associ\'ees \`a
une
forme modulaire et \`a un corps quadratique imaginaire, {\it J. Lon. Math.
Soc.}
{\bf 38} (1988) 1-32.
\smallskip

\noindent [PSR] Piatetski-Shapiro, I. I. and S. Rallis: Rankin triple
$L$-functions,
{\it Compositio Math.},  {\bf 64} (1987) 31-115.
\smallskip

\noindent [Sha]  Shahidi, F.: On the Ramanujan conjecture and finiteness
of poles for certain $L$-functions, {\it Annals of Math.}, {\bf 127} (1988)
547-584.
\smallskip

\noindent [S]  Shimura, G.: On a class of nearly holomorphic automorphic
forms, {\it Annals of Math.} {\bf 123}  (1986)  347-406.
\smallskip

\noindent [So]  Sofer, A.: $p$-adic interpolation of half-integral
weight modular forms, {\it Contemp. Math.}, {\bf 174}  (1995)  119-128.
\smallskip

\noindent [St]  Stevens, G.: $\Lambda$-adic modular forms of half-integral
weight and a $\Lambda$-adic Shintani lifting, {\it Contemp. Math.}, {\bf 174}
(1995)  129-151.
\smallskip

\noindent [V] Rodriguez Villegas, F.:  On the Taylor coefficients of theta
functions of CM elliptic curves, {\it Contemp. Math.}, {\bf 174} (1995)
185-201.

\smallskip

\noindent [W]  Waldspurger, J.-L.: Sur les valeurs de certaines fonctions
$L$
automorphes en leur centre de sym\'etrie, {\it Compositio Math.},
{\bf 54} (1985) 173-242.
\smallskip

\end

%$$D_H(\bold{f},g,h)(P) = P(\theta(\bold{f},g,h)) = P[ <h\cdot
%d\mu_g)),(T_{\bold{f}}\otimes 1)>] $$

$$\multline
D_H(\bold{f},g,h)(P)= H(P)\cdot 1_{\bold{f}_P}\cdot [h \cdot \int_{Z} x^r
\cdot d\mu_g] =
  H(P)\cdot 1_{\bold{f}_P}\cdot [h \cdot d^r(g_p)\tag 1.5.2.1
$$
by (1.4.2.2).